\documentclass{daj}

\usepackage{amsfonts}
\usepackage{amssymb}
\usepackage{amsthm}
\usepackage{amsmath}
\usepackage[capitalise]{cleveref}

\newcommand{\R}{\mathbb{R}}
\newcommand{\Z}{\mathbb{Z}}
\newcommand{\Span}{\text{\textup{span}}}
\newcommand*{\vol}{\mathop{\textup{vol}}\nolimits}

\newtheorem{prop}{Proposition}[section]
\newtheorem{theorem}[prop]{Theorem}
\newtheorem*{mink.2.thm}{Minkowski's second theorem}
\newtheorem{lemma}[prop]{Lemma}
\newtheorem{corollary}[prop]{Corollary}
\newtheorem{conjecture}[prop]{Conjecture}

\theoremstyle{definition}

\theoremstyle{remark}

\newtheorem*{remark*}{Remark}
\newtheorem{remark}[prop]{Remark}
\newtheorem*{remarks*}{Remarks}

\dajAUTHORdetails{%
  title = {New Bounds in the Discrete Analogue of Minkowski's Second Theorem}, 
  author = {Matthew Tointon},
  plaintextauthor = {Matthew Tointon},
}   

\dajEDITORdetails{%
   year={2024},
   number={7},
   received={23 March 2023},   
   published={1 October 2024},  
   doi={10.19086/da.122972},       
}   

\begin{document}

\begin{frontmatter}[classification=text]


\author[toint]{Matthew Tointon}

\begin{abstract}
We adapt an argument of Tao and Vu to show that if $\lambda_1\le\cdots\le\lambda_d$ are the successive minima of an origin-symmetric convex body $K$ with respect to some lattice $\Lambda<\R^d$, and if we set $k=\max\{j:\lambda_j\le1\}$, then $K$ contains at most $2^k(1+\frac{\lambda_k}2)^k/\lambda_1\cdots\lambda_k$ lattice points. This provides improved bounds in a conjecture of Betke, Henk and Wills (1993), and verifies that conjecture asymptotically as $\lambda_k\to0$. We also obtain a similar result without the symmetry assumption.
\end{abstract}
\end{frontmatter}


\section{Introduction}
An important and fundamental result in the geometry of numbers is Minkowski's famous \emph{second theorem}, from 1896, which bounds the volume of a convex body in Euclidean space in terms of the successive minima of that body with respect to a given lattice. In 1993, Betke, Henk and Wills conjectured a certain discrete analogue of this statement, with the volume of the convex body replaced by the number of lattice points it contains. The purpose of this note is to present some new bounds in the direction of this conjecture.

We start by recalling the necessary definitions. A \emph{convex body} $K\subseteq\R^d$ is a bounded convex set with non-empty interior. We call $K$ \emph{origin symmetric} if $K=-K$. Given a lattice $\Lambda<\R^d$, we define the \emph{successive minima} $\lambda_1\le\cdots\le\lambda_d$ of $K$ with respect to $\Lambda$ via
\[
\lambda_i=\inf\{\lambda>0:\dim\Span_\R(\textstyle\frac{\lambda}2(K-K)\cap\Lambda)\ge i\}.
\]
Note that if $K$ is origin symmetric then $K=\frac12(K-K)$, so that
\[
\lambda_i=\inf\{\lambda>0:\dim\Span_\R(\lambda K\cap\Lambda)\ge i\}.
\]
\begin{mink.2.thm}
Suppose that $K\subseteq\R^d$ is a convex body, $\Lambda<\R^d$ is a lattice, and $\lambda_{1}\le\cdots\le\lambda_{d}$ are the successive minima of $K$ with respect to $\Lambda$. Then
\[
\vol(K)\le\frac{2^d\det\Lambda}{\lambda_1\cdots\lambda_d}.
\]
\end{mink.2.thm}
Minkowski originally proved his theorem with the additional assumption that $K$ was origin symmetric, but the general case easily reduces to that case via the Brunn--Minkowski inequality, which implies in particular that $\vol(K)\le\vol(\frac12(K-K))$.

Minkowski's second theorem is very useful -- for instance, I myself first encountered it due to its application in proofs of Freiman's theorem on sets of integers with small sumset (see e.g. \cite[Theorem 3.1]{ruzsa}) -- and has attracted great interest, for example being given a number of different proofs (see e.g. \cite{henk} and the references therein).

Betke, Henk and Wills \cite[Conjecture 2.1]{bhw} conjectured an analogue of Minkowski's second theorem with the volume of $K$ replaced by the number of lattice points it contains. They originally stated their conjecture for origin-symmetric convex bodies only, but Malikiosis \cite{malikiosis} later conjectured that it should also hold without the symmetry assumption.
\begin{conjecture}[Betke, Henk and Wills; Malikiosis]\label{conj}
Suppose $K\subseteq\R^d$ is a convex body, $\Lambda<\R^d$ is a lattice, and $\lambda_{1}\le\cdots\le\lambda_{d}$ are the successive minima of $K$ with respect to $\Lambda$. Then
\[
\#(K\cap\Lambda)\le\prod_{i=1}^d\left\lfloor\frac2{\lambda_i}+1\right\rfloor.
\]
\end{conjecture}
Note that equality is achieved for $K=[0,a_1]\times\cdots\times[0,a_d]$ with $a_i>0$.

Actually, \cref{conj} is more than a mere analogue of Minkowski's second theorem: as Betke, Henk and Wills themselves noted in \cite[Proposition 2.2]{bhw}, one can easily recover Minkowski's second theorem from it by taking the limit
\begin{equation}\label{eq:mink2}
\frac{\vol(K)}{\det\Lambda}=\lim_{r\to0}r^d\#(K\cap r\Lambda)\le\lim_{r\to0}r^d\prod_{i=1}^d\left\lfloor\frac2{r\lambda_i}+1\right\rfloor=\frac{2^d}{\lambda_1\cdots\lambda_d}.
\end{equation}

\cref{conj} is trivial in dimension $1$. Betke, Henk and Wills proved it for origin-symmetric convex bodies in dimension $2$, and in general dimension $d$ showed that it holds up to a factor of roughly $d!$. Henk \cite{henk} later improved this to a factor of $2^{d-1}$. Malikiosis \cite{malikiosis} then proved the conjecture for arbitrary convex bodies in dimension at most $3$, and in general dimension $d$ showed that
\begin{equation}\label{eq:malik}
\#(K\cap\Lambda)\le\frac4e\left(\sqrt{3}\right)^{d-1}\prod_{i=1}^d\left\lfloor\frac2{\lambda_i}+1\right\rfloor,
\end{equation}
and even
\begin{equation}\label{eq:malik.sym}
\#(K\cap\Lambda)\le\frac4e\left(\sqrt[3]{\frac{40}9}\right)^{d-1}\prod_{i=1}^d\left\lfloor\frac2{\lambda_i}+1\right\rfloor
\end{equation}
if $K$ is assumed origin symmetric \cite[Theorem 3.1.2]{malikiosis}. Most recently, Freyer and Lucas \cite{freyer-lucas} proved that
\begin{equation}\label{eq:fl}
\#(K\cap\Lambda)\le\prod_{i=1}^d\left(\frac2{\lambda_i}+d\right),
\end{equation}
which improves on \eqref{eq:malik} in fixed dimension for small enough successive minima, and is the first bound in \cref{conj} strong enough to imply Minkowski's second theorem via the argument of \eqref{eq:mink2}. On the other hand, \eqref{eq:malik} remains stronger than \eqref{eq:fl} for example if the dimension goes to infinity with the successive minima bounded below.

In the present note we adapt a proof of Minkowski's second theorem due to Tao and Vu \cite[Theorem 3.30]{tao-vu} to obtain improvements to both \eqref{eq:malik} and \eqref{eq:malik.sym} in all cases, and to \eqref{eq:fl} in some cases, such as when those successive minima that are less than $2$ are all within a factor of $d$ of one another. Moreover, like Freyer and Lucas's bound \eqref{eq:fl}, these improvements are strong enough to imply Minkowski's second theorem via the argument of \eqref{eq:mink2}.

In the origin-symmetric case we have the following result.
\begin{theorem}\label{thm}
Suppose $K\subseteq\R^d$ is an origin-symmetric convex body, and $\Lambda<\R^d$ is a lattice. Let $\lambda_{1}\le\cdots\le\lambda_{d}$ be the successive minima of $K$ with respect to $\Lambda$. If $\lambda_1>1$ then $\#(K\cap\Lambda)=1$. Otherwise, setting $k=\max\{j:\lambda_j\le1\}$, for any real numbers $\mu_1\le\cdots\le\mu_k$ satisfying $0<\mu_i\le\lambda_i$ we have
\[
\#(K\cap\Lambda)\le\frac{2^k(1+\frac{\mu_k}2)^k}{\mu_1\cdots\mu_k}.
\]
\end{theorem}
The statement given in the abstract is what one obtains from this result by setting each $\mu_i=\lambda_i$. This already improves on \eqref{eq:malik.sym}, but depending on the values of the successive minima, the flexibility to choose the $\mu_i$ can improve the bound still further, particularly in high dimension if there is one successive minimum that is larger than a large number of the other ones. For example, if $\lambda_d=1$ and all other $\lambda_i$ are at most $\frac12$, then taking $\mu_d=\lambda_d=1$ results in a bound of $3^d/\lambda_1\cdots\lambda_{d-1}$, whereas taking $\mu_i=\lambda_i$ for $i<d$ and $\mu_d=\frac12$ gives a bound of $2(\frac52)^d/\lambda_1\cdots\lambda_{d-1}$. As we will see in \cref{cor:asym} below, in the asymmetric setting this flexibility is crucial for improving on \eqref{eq:malik}.

The assertion that $\#(K\cap\Lambda)=1$ if $\lambda_1>1$ is of course trivial, but separating that case makes the definition of $k$ less cumbersome.

In the general setting we have the following result.
\begin{theorem}\label{thm:asym}
Suppose $K\subseteq\R^d$ is a convex body, and $\Lambda<\R^d$ is a lattice. Let $\lambda_{1}\le\cdots\le\lambda_{d}$ be the successive minima of $K$ with respect to $\Lambda$. If $\lambda_1>2$ then $\#(K\cap\Lambda)\le1$. Otherwise, setting $k=\max\{j:\lambda_j\le2\}$, for any real numbers $\mu_1\le\cdots\le\mu_k$ satisfying $0<\mu_i\le\lambda_i$ we have
\[
\#(K\cap\Lambda)\le\frac{2^k(1+\frac{\mu_k}2)^k}{\mu_1\cdots\mu_k}.
\]
\end{theorem}
Unlike in the symmetric setting, an asymmetric convex body can contain several lattice points even if its successive minima are all greater than $1$, as is the case with $[0,1]^d\subseteq\R^d$ and the lattice $\Z^d$, for example. Nonetheless, if $\lambda_1>2$ then $(K-K)\cap\Lambda=\{0\}$, and so $\#(K\cap\Lambda)\le1$ as claimed. This is also essentially the reason for the $1$ in the definition of $k$ in \cref{thm} being replaced by $2$ in \cref{thm:asym}.

Setting $\mu_i=\lambda_i$ in \cref{thm:asym} immediately improves on \eqref{eq:malik} if there is no $i$ such that $1<\lambda_i\le2$ (and in various other cases). However, \eqref{eq:malik} may still be stronger if there exists such an $i$. For example, if $\lambda_1=\cdots=\lambda_{d-1}=1/M$ for some $M>1$ and $\lambda_d=2$ then setting $\mu_i=\lambda_i$ gives a bound of $\frac124^dM^{d-1}$, whilst \eqref{eq:malik} gives a bound of roughly $\frac8e(2\sqrt{3}M)^{d-1}$. Nonetheless, by choosing $\mu_i=1$ whenever $1\le\lambda_i\le2$ we obtain the following corollary.
\begin{corollary}\label{cor:asym}
Suppose $K\subseteq\R^d$ is a convex body and $\Lambda<\R^d$ is a lattice. Let $\lambda_{1}\le\cdots\le\lambda_{d}$ be the successive minima of $K$ with respect to $\Lambda$, and suppose there exists $\lambda_i\in(1,2]$. Then, setting $k=\max\{j:\lambda_j\le1\}$ and $m=\max\{j:\lambda_j\le2\}$, we have
\[
\#(K\cap\Lambda)\le\frac{3^m}{\lambda_1\cdots\lambda_k}.
\]
\end{corollary}
To see that \cref{cor:asym} improves on \eqref{eq:malik} in dimension $2$ or greater whenever it applies, write the bound of \eqref{eq:malik} as $A_1\cdots A_d$, with $A_1=\frac4e\lfloor\frac2{\lambda_1}+1\rfloor$ and $A_i=\sqrt{3}\lfloor\frac2{\lambda_i}+1\rfloor$ for $i>1$, and write the bound of \cref{cor:asym} as $B_1\cdots B_m$, with $B_i=3/\lambda_i$ for $i\le k$ and $B_{k+1},\ldots,B_m=3$. Certainly $A_i>1$ for $i>m$. Also, $B_m=3$ and $B_1=3\max\{1,\frac1{\lambda_1}\}$, whilst $A_m=2\sqrt{3}$ and $A_1\ge\frac8e\max\{1,\frac1{\lambda_1}\}$, so that $B_1B_m<A_1A_m$. Finally, for all other $i$ we have $B_i=3\max\{1,\frac1{\lambda_i}\}$ and $A_i\ge2\sqrt{3}\max\{1,\frac1{\lambda_i}\}$, and hence $B_i<A_i$.

In general, the best values of $\mu_i$ to take in \cref{thm,thm:asym} will depend on the exact values of the successive minima. In particular, \cref{cor:asym} is not meant to represent an optimal choice of $\mu_i$ in \cref{thm:asym}, but rather just one choice that demonstrates that \cref{thm:asym} can always improve on \eqref{eq:malik}.

Finally, before embarking on the proofs, it is perhaps worth remarking that \cref{conj} can be seen as part of a broader movement to discretise inequalities from convex geometry; see \cite{MR4297379}, for example.

\section{Proof of the theorems}

The main technical ingredient in Tao and Vu's proof of Minkowski's second theorem, and one which plays a central role in the proofs of \cref{thm,thm:asym}, is a lemma they call the \emph{squeezing lemma}.
\begin{lemma}[squeezing lemma {\cite[Lemma 3.31]{tao-vu}}]\label{lem:squeezing}
Suppose $K\subseteq\R^d$ is an open convex body, $V\le\R^d$ is a $j$-dimensional subspace, and $A\subseteq K$ is an open set. Then there exists an open subset $A'\subseteq K$ satisfying $\vol(A')=\mu^j\vol(A)$ and $(A'-A')\cap V\subseteq(\mu(A-A))\cap V$.
\end{lemma}
The statement of \cite[Lemma 3.31]{tao-vu} actually includes an assumption that $K$ is origin-symmetric, but this is not used anywhere in the proof. A reader seeking reassurance that this is the case can also consult \cite[Lemma 3.5.2]{tointon.book}, where the lemma is stated exactly as above, and proved following Tao and Vu.

In the following proof, we denote by $K^\circ$ the interior of $K$ and by $\overline K$ the closure of $K$.

\begin{proof}[Proof of \cref{thm,thm:asym}]
Since the proofs of \cref{thm,thm:asym} have considerable overlap, we prove these theorems simultaneously. Thus, throughout this proof, $K$ is a convex body, not necessarily origin symmetric except in the parts of the proof that apply specifically to \cref{thm}.

We may assume that $\#(K\cap\Lambda)\ge2$. On translating by an element of $\Lambda$, we may therefore also assume that $0\in K$, and that $K\cap\Lambda$ generates a non-zero subspace $V$ of $\R^d$, say of dimension $m$. Let $\lambda_1'\le\cdots\le\lambda_m'$ be the successive minima of the convex body $K'=K\cap V$ with respect to the lattice $\Lambda'=\Lambda\cap V$ in $V$, noting that $\lambda_i'\ge\lambda_i$ for $i=1,\ldots,m$. The fact that $0\in K$ implies that $(K'-K')\cap\Lambda'\supseteq K'\cap\Lambda$, so that $(K'-K')\cap\Lambda$ spans $V$ and $\lambda_m'\le2$. Moreover, if $K$ is origin symmetric then so is $K'$, so that $\frac12(K'-K')\cap\Lambda'=K'\cap\Lambda'$ spans $V$ and $\lambda_m\le1$. On considering $K'$ and $\Lambda'$ in place of $K$ and $\Lambda$, we may therefore assume that $k=d$.

Choose inductively a list $v_1,\ldots,v_d\in\Lambda$ of linearly independent vectors such that $v_1,\ldots,v_i\in\frac{\lambda_i}2(\overline K-\overline K)$ for every $i$. For each $i=0,1,\ldots,d$ set $V_{i}=\Span_\R(v_1,\ldots,v_i)$ and $\Lambda_{i}=\Lambda\cap(V_{i}\setminus V_{i-1})$, noting that $\frac{\mu_j}2(K^\circ-K^\circ)\cap\Lambda_j=\{0\}$ for every $\mu_j\in(0,\lambda_j]$.

Following Tao and Vu \cite[\S3.5]{tao-vu} (but with $\mu_i$ in place of $\lambda_i$), set $K_0=\frac{\mu_d}{2}K^\circ$ and then, starting with $A_0=K_0$, apply \cref{lem:squeezing} iteratively to obtain a sequence $A_0,A_1,\ldots,A_{d-1}$ of open subsets of the open convex body $K_0$ such that
\[
\vol(A_i)=\left(\frac{\mu_i}{\mu_{i+1}}\right)^i\vol(A_{i-1})
\]
and
\[
(A_i-A_i)\cap V_i\subseteq\left(\frac{\mu_i}{\mu_{i+1}}(A_{i-1}-A_{i-1})\right)\cap V_i
\]
for each $i$. Note that
\[
\vol(A_{d-1})=\frac{\mu_{1}\ldots\mu_{d}}{2^d}\vol(K).
\]
Furthermore, for each $j$ we have
\begin{align*}
(A_{d-1}-A_{d-1})\cap V_j&\subseteq\textstyle\frac{\mu_j}{\mu_d}(A_{j-1}-A_{j-1})\cap V_j\\
     &\subseteq\textstyle\frac{\mu_j}{\mu_d}(K_0-K_0)\cap V_j\\
    &=\textstyle\frac{\mu_j}2\left(K^\circ-K^\circ\right)\cap V_j,
\end{align*}
and hence $(A_{d-1}-A_{d-1})\cap\Lambda_j\subseteq\frac{\mu_j}2(K^\circ-K^\circ)\cap\Lambda_j=\{0\}$. Since this holds for all $j$, we conclude that $(A_{d-1}-A_{d-1})\cap\Lambda=\{0\}$.

Departing now from Tao and Vu's argument, note that this last property implies that $\vol(X+A_{d-1})=\vol(A_{d-1})\#X$ for any finite subset $X\subseteq\Lambda$, and in particular that
\[
\vol((K\cap\Lambda)+A_{d-1})=\vol(A_{d-1})\#(K\cap\Lambda)=\frac{\mu_1\cdots\mu_d}{2^d}\vol(K)\#(K\cap\Lambda).
\]
Since
\[
(K\cap\Lambda)+A_{d-1}\subseteq K+K_0\subseteq\left(1+\frac{\mu_d}2\right)K,
\]
it follows that
\[
\#(K\cap\Lambda)\le\frac{2^d(1+\frac{\mu_d}2)^d}{\mu_1\cdots\mu_d},
\]
as required.
\end{proof}

\begin{remark}
The same reduction we performed in the second paragraph of the above proof shows that one may replace $d$ by $k$ in each of \eqref{eq:malik}, \eqref{eq:malik.sym} and \eqref{eq:fl}.
\end{remark}

\begin{remark*}
Although Minkowski's second theorem follows from \cref{thm} via the argument of \eqref{eq:mink2}, it follows more directly from applying Blichfeldt's lemma to the set $A_{d-1}$ with $\mu_i=\lambda_i$, as Tao and Vu do. Note that, in the asymmetric case, this is also `direct' in the sense that it does not necessitate passing via the symmetric case using the Brunn--Minkowski inequality.
\end{remark*}

\section*{Acknowledgments} 
I thank Martin Henk for the suggestion to consider asymmetric convex bodies, which significantly improved this note, and for help with the references. I also thank both him and an anonymous referee for comments on earlier versions of the manuscript.

\bibliographystyle{amsplain}


\begin{dajauthors}
\begin{authorinfo}[toint]
  Matthew Tointon\\
  School of Mathematics\\
  University of Bristol\\
  United Kingdom\\
  m\imagedot{}tointon\imageat{}bristol\imagedot{}ac\imagedot{}uk \\
  \url{https://tointon.neocities.org/}
\end{authorinfo}
\end{dajauthors}

\end{document}